\documentclass[11pt,a4paper]{article}

\usepackage{epsfig,amsthm,amscd,amsxtra,amsbsy}
\usepackage{amsmath,amssymb,latexsym, amsfonts}

\begin{document}

\newcommand{\bp}{{\mbox{}\hskip 0.1cm}}
\newcommand{\ba}{\begin{equation}}
\newcommand{\ea}{\end{equation}}

\newtheorem{theorem}{Theorem}[section]
\newtheorem{lemma}{Lemma}[section]
\newtheorem{proposition}{Proposition}[section]
\newtheorem{definition}{Definition}[section]
\newtheorem{assumption}{Assumption}[section]

\textwidth 180mm
\textheight 210mm

\baselineskip 0.6cm
\parskip 0.0cm
\parindent 0.7cm

\pagenumbering{arabic}

\title{\bp {Multi-symplectic Birkhoffian Structure for PDEs with Dissipation Terms }\thanks{
\quad Supported by the Special Funds for Major State Basic
Research Projects , G 1999, 032800}}
\author{Hongling Su\quad\quad\quad
 Mengzhao Qin\\
 CAST (World Laboratory), \\
Institute of Computational Mathematics\\
 and Scientific/Engineering Computing, \\
Academy of Mathematics and System Sciences,\\
Chinese Academy of Sciences, Beijing, 100080, China.}
\maketitle
\begin{abstract}The multi-symplectic form for Hamiltonian PDEs leads to a general framework for geometric numerical
schemes that preserve a discrete version of the conservation of
symplecticity. The cases for systems or PDEs with dissipation
terms has never been extended. In this paper, we suggest a new
extension for generalizing the multi-symplectic form for
Hamiltonian systems to systems with dissipation which never have
remarkable energy and momentum conservation properties. The
central idea is that the PDEs is of a first-order type that has a
symplectic structure depended explicitly on time variable, and
decomposed into distinct components representing space and time
directions. This suggest a natural definition of multi-symplectic
Birkhoff's equation as a multi-symplectic structure from that a
multi-symplectic dissipation law is constructed. We  show that
this definition leads to deeper understanding  relationship
between functional principle and PDEs. The concept of
multi-symplectic integrator is also discussed.
\end{abstract}
\vskip 0.3cm
\setcounter{page}{1}
\section{\Large\bf Introduction}
The multi-symplectic Hamilton's Equation is proposed first time by
Bridges in 1997\cite{Bridges1}. Birkhoff's Equation is named by
Santilli\cite{Santilli2}, who just discussed $1$-dimensional
cases. We return to  multi-symplectic Hamilton's Equation, and we
just discuss $1+1$-dimensional cases. It is a local concept of
symplecticity for Hamiltonian systems. The concept can be
abstracted as follows. Let $M$ and $K$ be any skew symmetric
matrices on $\mathbb{R}^n (n\geq 3)$ and Let $S:\quad
\mathbb{R}^n\longrightarrow \mathbb{R}$ be any smooth function.
Then we will call a system of the following form a Hamiltonian
system on  a multi-symplectic structure, \ba\label{ds1.1}
Mz_t+Kz_x=\bigtriangledown_zS(z,x,t),\qquad z\in\mathbb{R}^n,
\quad(x,t)\in\mathbb{R}^2, \ea where the general gradient
$\bigtriangledown_z$ is defined with respect to the standard inner
product on $\mathbb{R}^n$, denoted by $<\cdot,\cdot>$. Associated
$M$ and $K$ are the so-called pre-symplectic forms
\ba\label{ds1.2} \omega(U,V)=<MU,V>,\qquad\kappa(U,V)=<KU,V>. \ea
The two $2$-forms $\omega$ and $\kappa$ define two symplectic
structure on $\mathbb{R}^m\quad(m=\text{rank}M\leq n)$ and
$\mathbb{R}^k\quad(k=\text{rank}K\leq n)$ respectively. The pair
$(\omega,\kappa)$ satisfies conservation of symplecticity
\ba\label{ds1.3} \frac{\partial}{\partial
t}\omega+\frac{\partial}{\partial x}\kappa=0. \ea A
multi-symplectic integrator is numerical approximation of equation
($\ref{ds1.1}$) \ba\label{ds1.4} M\partial
t^{i,j}z_{i,j}+K\partial
x^{i,j}z_{i,j}=(\bigtriangledown_zS(z_{i,j},x_i,t_j))_{i,j}, \quad
\text{where}\quad z_{i,j}=z(x_i,t_j), \ea which has a discrete
conservation law as \ba\label{ds1.5}
\partial t^{i,j}\omega_{i,j}+\partial x^{i,j}\kappa_{i,j}=0.
\ea
\par
The problem with this type to abstract and to discrete the systems
or PDEs is that it is subject to conservative systems or
nonconservative systems without dissipation terms. A more general
multi-symplectic structure is to represent some nonconservative
systems with dissipation terms, since,---different from form
($\ref{ds1.1}$)---,symplecticity may depend on time variables
explicitly, and the right hand is a gradient function (vector
function) plus a dissipation terms.
\par
Similarly, let $M(z,x,t)$ and $K(z,x,t)$ are skew symmetric
metrics on $\mathbb{R}^n\geq 2$, and $B;\quad
\mathbb{R}^n\times\mathbb{R}^2\longrightarrow\mathbb{R}$ be any
smooth function. Then a Birkhoffian system under the following
multi-symplectic structure is given \ba\label{ds1.6}
M(z,x,t)z_t+K(z,x,t)z_x=\bigtriangledown_z B(z,x,t)+\frac{\partial
F}{\partial t}+\frac{\partial G}{\partial x}, \ea where $M$ and
$F$, $K$ and $G$ are independent. We will show their relationship
in the next section.
\par
Structure ($\ref{ds1.6}$) is interesting for two reasons: it is a
natural generalization of structure ($\ref{ds1.1}$) as well as
Birkhoff's representation for ODEs to PDEs, and it can be used to
develop symplectic numerical integrators for nonconservative
systems with dissipation terms. For example, consider the equation
describe a linear damped string \ba\label{ds1.7}
u_{tt}-u_{xx}+u+\alpha u_t+\beta u_x=0. \ea we introduce the new
variables $p$ and $q$, defined below, and equation ($\ref{ds1.7}$)
is rewritten as first-order form \ba\label{ds1.8} u_t=p,\qquad
u_x=q,\quad p_t-q_x+u+2p=0, \ea or, with
$z=(u,p,q)^{\text{T}}\in\mathbb{R}^3$, we obtain a
multi-symplectic Birkhoff's equations (ds1.6), where
\ba\label{ds1.9} M=\begin{pmatrix}
{0}&\text{e}^{\alpha t-\beta x}&0\\
{-\text{e}^{\alpha t-\beta x}}&0&0\\
{0}&0&0
\end{pmatrix},\qquad\qquad
K=\begin{pmatrix}
{0}&0&-\text{e}^{\alpha t-\beta x}\\
{0}&0&0\\
{\text{e}^{\alpha t-\beta x}}&0&0
\end{pmatrix},
\ea
and
\ba\label{ds1.10}
\begin{split}
F=(-\frac{1}{2}e^{\alpha t-\beta x}p,\frac{1}{2}e^{\alpha t-\beta x}u,0)^{\text{T}},0),\\
G=(\frac{1}{2}e^{\alpha t-\beta x}q,0,\frac{1}{2}e^{\alpha t-\beta x}u,0)^{\text{T}}),\\
B=-\frac{1}{2}\text{e}^{\alpha t-\beta x}[u^2+p^2-q^2+\alpha up+\beta uq].
\end{split}
\ea
\par
We now show an outline of the paper. In section 2 we introduce
accurate definition of multi-symplectic Birkhoff's equations. From
the point of inverse variational problem we discuss the conditions
of self-adjointness for Birkhoff's equations. This leads naturally
to a variational principle. In section 3, we introduce continuous
and discrete dissipation laws of symplecticity, energy and
momentum, then give definition of multi-symplectic Birkhoffian
integrator. In section 4,we consider equation ($\ref{ds1.7}$) to
illustrate how a multi-symplectic Birkhoffian integrator to
preserve the discrete version of the  dissipation laws.
\section{Multi-symplectic Birkhoff's Equation}
The general multi-symplectic structure we shall study is given by
\ba\label{ds2.1}
\begin{split}
&[\frac{\partial F_\nu(z,x,t)}{\partial z^\mu}-\frac{\partial F_\mu(z,x,t)}{\partial z^\nu}]\frac{\partial z_\nu}{\partial t}+[\frac{\partial G_\nu(z,x,t)}{\partial z^\mu}-\frac{\partial G_\mu(z,x,t)}{\partial z^\nu}]\frac{\partial z_\nu}{\partial x}\\
&-[\frac{\partial B(z,x,t)}{\partial z^\mu}+\frac{\partial F_\mu(z,x,t)}{\partial t}+\frac{\partial G_\mu(z,x,t)}{\partial x}]=0.
\end{split}
\ea
Equations ($\ref{ds2.1}$) are called multi-symplectic Birkhoff's equations. Why using this name is reviewed in Ref.\cite{Santilli2}. The function $B(z,x,t)$ is called Birkhoffian. A representation of a system via multi-symplectic Birkhoff's equations is called a multi-symplectic Birkhoffian representation, and the system is called multi-symplectic system.
\par
Now we consider a general first-order form in the self-explanatory unified notation
\ba\label{ds2.2}
M_{\mu\nu}(z,x,t)z^\nu_t+K_{\mu\nu}z^\nu_x+D_mu(z,x,t)=0.
\ea
We give the conditions of variational self-adjointness for equation ($\ref{ds2.2}$) as follows.
\begin{theorem}\label{dsth2.1}
(self-adjointness of the covariant general form) A necessary and sufficient condition for equations ($\ref{ds2.1}$) to be self-adjoint in a region $\widetilde{\Re}$ of points (z,x,t) is that all the following conditions
\ba\label{dsth2.1.1}
\begin{split}
M_{\mu\nu}+M_{\nu\mu}=0,\qquad K_{\mu\nu}+K_{\nu\mu}=0, \\
\frac{\partial M_{\mu\nu}}{\partial z^\tau}+\frac{\partial M_{\nu\tau}}{\partial z^\mu}+\frac{\partial M_{\tau\mu}}{\partial z^\nu},\qquad \frac{\partial K_{\mu\nu}}{\partial z^\tau}+\frac{\partial K_{\nu\tau}}{\partial z^\mu}+\frac{\partial K_{\tau\mu}}{\partial z^\nu}, \\
\frac{\partial M_{\mu\nu}}{\partial t}=\frac{\partial D}{\partial z^\nu}-\frac{\partial D}{\partial z^\mu},\qquad \frac{\partial K_{\mu\nu}}{\partial t}=\frac{\partial D}{\partial z^\nu}-\frac{\partial D}{\partial z^\mu},
\end{split}
\ea
are identically verified in $\widetilde{\Re}$.
\end{theorem}
\par
Note that the conditions of self-adjointness are formulated for the covariant form because this is the form derivable from a variational principle which is the most general possible first-order action functional given by the pfaffian action
\ba\label{ds2.3}
\pounds(\widetilde{E})=\int\limits^{t_2}_{t_1}\int\limits^{x_2}_{x_1}\text{d}t\text{d}x[F\nu(z,x,t)z^\nu_t+G_\nu(z,x,t)z^\nu_x-B(z,x,t)](\widetilde{E}).
\ea

Strip away the integral from $\pounds(\widetilde{E})$ in ($\ref{ds2.3}$), define a general Lagrange density $\mathcal{L}=F\nu(z,x,t)z^\nu_t+G_\nu(z,x,t)z^\nu_x-B(z,x,t)$. It doesn't exist a Legendre transformation (including conjugate momenta ) leads reduction from the Lagrange equation to Hamilton's equation.  While from the variational principle we can write a Birkhoff's equation directly, and for any self-adjoint covariant system, a variational principle can be constructed.
\begin{proposition}\label{dsth2.2}
(Self-adjointness of Birkhoff's Equations). Necessary and sufficient condition for a general first order system ($ds2.2$) to be self-adjoint is that it is of the Birkhoffian type, i.e.,
\ba\label{dsth2.2.1}
\begin{split}
&M_{\mu\nu}(z,x,t)z^\nu_t+K_{\mu\nu}z^\nu_x+D_mu(z,x,t)=\\
&[\frac{\partial F_\nu(z,x,t)}{\partial z^\mu}-\frac{\partial F_\mu(z,x,t)}{\partial z^\nu}]\frac{\partial z_\nu}{\partial t}+[\frac{\partial G_\nu(z,x,t)}{\partial z^\mu}-\frac{\partial G_\mu(z,x,t)}{\partial z^\nu}]\frac{\partial z_\nu}{\partial x}\\
&-[\frac{\partial B(z,x,t)}{\partial z^\mu}+\frac{\partial F_\mu(z,x,t)}{\partial t}+\frac{\partial G_\mu(z,x,t)}{\partial x}]=0.
\end{split}
\ea
\end{proposition}

\section{Multi-symplecticity Dissipation Laws and Multi-symplectic Integrators}
The multi-symplectic Birkhoff's equations ($\ref{ds2.1}$) determines 2 $two$-forms $\omega$ and $kappa$. Since
\ba\label{ds3.1}
\begin{split}
\frac{\partial}{\partial t}\omega(U,V)=<MU_t,V>+<MU,V_t>+<M_tU,V>,\\
\frac{\partial}{\partial x}\kappa(U,V)=<KU_x,V>+<KU,V_x>+<K_xU,V>,
\end{split}
\ea
one obtains
\ba\label{ds3.2}
\begin{split}
\frac{\partial}{\partial t}\omega(U,V)+\frac{\partial}{\partial x}\kappa(U,V)&=<MU_t+KU_x,V>-<U,MV_t+KV_x>+<M_tU+K_xU,V>\\
&=<F_{zt}U+G_{zx}U,V>-<F_{zt}V+G_{zx}V,U>\\
&=0.
\end{split}
\ea
This statement is true for any solution pair $(U,V)$ and we write is as
\ba\label{ds3.3}
\partial_t\omega+\partial_x\kappa=0.
\ea
 This is different from the conservation law ($\ref{ds1.3}$), ---although they look like each other---, since that $M_t\ne0$ and $K_x\ne0$. Therefore we define this law as a multi-symplectic dissipation law for Birkhoffian system, and in the next section this definition will be verified very natural.
\par
The numerical framework is start from the Birkhoff's equation ($\ref{ds2.1}$) with the associated dissipation law of symplecticity ($\ref{ds3.3}$). The idea is to develop a numerical scheme to approximate equations ($\ref{ds2.1}$) and it also preserve a discrete version of dissipation of symplecticity. A numerical discretization of ($\ref{ds3.3}$) can be written schematically as
\ba\label{ds3.4}
\begin{split}
M_{i,j}\partial_t^{i,j}z_{i,j}+K_{i,j}\partial_x^{i,j}z_{i,j}&=(\bigtriangledown_zS(z_{i,j},x_i,t_j))_{i,j}+(\partial_tF(z_{i,j},x_i,t_j))_{i,j}\\
&+(\partial_xG(z_{i,j},x_i,t_j))_{i,j},\qquad z_{i,j}=z(x_i,t_j),
\end{split}
\ea
and $\partial_t^{i,j}$ and $\partial_x^{i,j}$ are discretizations of the derivatives $\partial_t$ and $\partial_x$ respectively. A discrete version of the dissipation law can be written as follows
\ba\label{ds3.5}
\partial_t^{i,j}\omega_{i,j}+\partial_x^{i,j}\kappa_{i,j}=0,
\ea
where
\ba\label{ds3.6}
\omega_{i,j}=<M_{i,j}U_{i,j},V_{i,j}>,\qquad\kappa_{i,j}=<K_{i,j}U_{i,j},V_{i,j}>.
\ea
\begin{definition}
The numerical scheme determined by equation ($ref{ds3.5}$) is called a multi-symplectic integrator if discrete law ($\ref{ds3.6}$) is a dissipation law for equations ($\ref{ds3.5}$).
\end{definition}
\section{example}
We consider the following equation
\ba\label{ds4.1}
u_{tt}-u_{xx}+u+2u_t=0.
\ea
it has a multi-symplectic Birkhoffian representation that has structure
\ba\label{ds4.2}
M=\begin{pmatrix}
{0}&\text{e}^{2t}&0\\
{-\text{e}^{2t}}&0&0\\
{0}&0&0
\end{pmatrix},\qquad\qquad
K=\begin{pmatrix}
{0}&0&-\text{e}^{2t}\\
{0}&0&0\\
{\text{e}^{2t}}&0&0
\end{pmatrix},
\ea
and
\ba\label{ds4.3}
-\frac{1}{2}\text{e}^{2t}[u^2+p^2-q^2+2up].
\ea
It has a dissipation law
\ba\label{ds4.4}
\frac{\text{d}}{\text{d}t}(\text{d}p\wedge\text{d}u)+\frac{\text{d}}{\text{d}x}(\text{d}q\wedge\text{d}u)=-2\text{d}p\wedge\text{d}u,
\ea
using the conditions of self-adjointness, the dissipation law can be rewritten as\ba\label{ds4.5}
frac{\text{d}}{\text{d}t}(\text{e}^{2t}\text{d}p\wedge\text{d}u)+\frac{\text{d}}{\text{d}x}(\text{e}^{2t}\text{d}q\wedge\text{d}u)=0.
\ea

Approximating the equation ($\ref{ds4.1}$) in Birkhoffian structure as follows
\ba\label{ds4.6}
\begin{split}
\frac{\text{e}^{\delta t}p_{i+1/2,j+1}-\text{e}^{-\delta t}p_{i+1/2,j}}{\delta t}+\frac{q_{i+1,j+1/2}-q_{i,j+1/2}}{\delta x}=-u_{i+1/2,j+1/2},\\
\frac{u_{i+1/2,j+1}-u_{i+1/2,j}}{\delta t}=p_{i+1/2,j+1/2},\\
\frac{u_{i+1,j+1/2}-u_{i,j+1/2}}{\delta t}=q_{i+1/2,j+1/2},\\
\end{split}
\ea
where
\ba\label{ds4.7}
\begin{split}
p_{i+1/2,j+1/2}=\frac{1}{4}\text{e}^{\delta t}(p_{i+1,j+1}+p_{i,j+1})+\frac{1}{4}\text{e}^{-\delta t}(p_{i+1,j}+p_{i,j}),\\
u_{i+1/2,j+1/2}=\frac{1}{4}(u_{i+1,j+1}+u_{i,j+1}+u_{i+1,j}+u_{i,j}),\\
q_{i+1/2,j+1/2}=\frac{1}{4}(u_{i+1,j+1}+u_{i,j+1}+u_{i+1,j}+u_{i,j}).
\end{split}
\ea
Discrete the dissipation law ($\ref{ds4.2}$) directly by using the same scheme, we get
\ba\label{ds4.8}
\begin{split}
&\frac{\text{e}^{\delta t}\text{d}p_{i+1/2,j+1}\wedge\text{d}u_{i+1/2,j+1}-\text{e}^{-\delta t}\text{d}p_{i+1/2,j}\wedge\text{d}u_{i+1/2,j}}{\delta t}\\
&-\frac{\text{d}q_{i+1,j+1/2}\wedge\text{d}u_{i+1,j+1/2}-\text{d}q_{i,j+1/2}\wedge\text{d}u_{i,j+1/2}}{\delta x}=0.
\end{split}
\ea
From the discrete equations ($\ref{ds4.4}$), we also can construct the dissipation, i.e. ($\ref{ds4.8}$) is a symplectic dissipation law for equations ($\ref{ds4.4}$). So the scheme determined by equations ($\ref{ds4.4}$) is a multi-symplectic integrator for Birkhoff's equations of system ($\ref{ds4.1}$).

\end{document}